\documentclass[11pt]{article}
\usepackage{amssymb,graphics,amsmath}
\usepackage[all]{xy}
\newtheorem{proposition}{Proposition}[section]
\newtheorem{theorem}{Theorem}[section]
\newtheorem{lemma}{Lemma}[section]

\newtheorem{conjecture}{Conjecture}[section]

\def\Z{\mathbb{Z}}

\def\proof{{\bf Proof) }}
\begin{document}

\title{Lens spaces given from L-space homology 3-spheres}
\author{Motoo Tange\\
\vspace{10pt}\\
  {\small Department of Mathematics, Osaka University,}\\
  {\small Machikaneyama 1-1 Toyonaka, Osaka \textup{560-0043}, Japan}\\
    {\small e-mail\textup{: \texttt{tange@math.sci.osaka-u.ac.jp}}}}
\date{}
\maketitle 
\abstract
We consider the problem when lens spaces are given from homology spheres, and 
demonstrate that many lens spaces are obtained from L-space homology sphere 
which the correction term $d(Y)$ is equal to $2$.
We show an inequality of slope and genus when $Y$ is L-space and $Y_p(K)$ is lens space.
\section{Introduction}
In the present paper we define lens space $L(p,q)$ to be the $p/q$-Dehn surgery of unknot, where $p,q$ are 
coprime integers.
Note that this orientation is opposite to usual one.
Let $K\subset S^3$ be a knot whose Dehn surgery is homeomorphic to a lens space.
Then we call that $K$ admits {\it lens surgery on $S^3$} or simply {\it lens surgery}.
The main problems of lens surgery involve when a lens space is obtained from Dehn surgery 
of a classical knot or when a knot $K$ admits lens surgery.

The Alexander polynomial $\Delta_K(t)$ of classical knot $K$ is useful for lens surgery problem.
For example for a certain knot $K$ yielding a lens space if $\Delta_K(t)$ is trivial, then $K$ is unknot, since the degree of
Alexander polynomial coincides with the genus of $K$.
If $\Delta_K(t)$ is $t-1+t^{-1}$, then $K$ is the trefoil knot in \cite{[4]}.
Although it is unknown whether the uniqueness holds otherwise, it seems that Alexander polynomial 
is an effective invariant for lens surgery.
As showed in \cite{[5]}, the doubly primitive knots which yields $L(p,q)$ can be distinguished
by the Alexander polynomials.

Many people have been done research on condition for lens surgery.
In \cite{[2]} J. Berge has defined doubly primitive knots, each of which is a classical knot yielding lens space.
Conversely he conjectured that these knots are all knots admitting lens surgery.
He had fallen the doubly primitive knots into several types.
But unfortunately the classification is not sure whether it is complete.
If the uniqueness above holds, then his conjecture is solved.

P. Ozsv\'ath and Z. Szab\'o have proven the constraints of Alexander polynomial $\Delta_K(t)$ of knots yielding lens space
by using Heegaard Floer homology.
The formula in Section 10.3 or Theorem~7.2 in \cite{[3]} represents the relationship between correction term and 
Alexander polynomial.
This formula is due to the surgery exact triangle among $S^3$, $S^3_0(K)$ and $L(p,q)$.

In \cite{[1]} R. Fintushel and R. Stern have shown a criterion for obtaining a lens space from a Dehn surgery of 
a homology sphere.
The statement is the following.
\begin{theorem}[{\cite{[1]}}]
The lens space $L(p,q)$ can be obtained as integral positive surgery on a knot in an integral homology 
three-sphere if and only if $q$ is a square modulo $p$.
\end{theorem}
On the other hand P. Ozsv\'ath and Z. Szab\'o have shown that for positive integers $k$ not divisible by $4$, 
$L(2k(3+8k),2k+1)$ cannot be obtained from Dehn surgery of $S^3$, but can be obtained from Dehn 
surgery of a Brieskorn homology sphere and moreover the lens space is obtained from 2-component link.
The detailed proof is in \cite{[3]}.

We would like to consider when a lens space can be obtained from 
a Dehn surgery of a knot in a homology sphere.
This issue seems to be worth considering but it does not seem that such 
studies about Dehn surgery of homology spheres
have been developed remarkably.
This present paper attends to Dehn surgeries of L-space homology spheres, which is defined as the manifolds 
whose Heegaard Floer homology is isomorphic to that of $S^3$.
The Heegaard Floer theoretic arguments of Dehn surgery of L-space homology sphere are the same as 
$S^3$ up to the absolute grading.
They lead to the three results: the surgery formula between the correction term and the Alexander polynomial
(Equation~(\ref{correction})) , 
the constraint of Alexander polynomial of lens surgery knot (Lemma~\ref{degg}), and the 
lower bound of the slope by the genus of the knot.
The lower bound is
$$2g(K)-1\le p.$$
Here we state an upper bounds, where it generalizes the main theorem in \cite{[8]}.
The notation $d(Y)$ is the correction term of $Y$, whose definition is in \cite{[3]}.
\begin{theorem}
\label{main1}
Let $Y$ be an L-space homology sphere.
Suppose that $Y_p(K)$ is a lens space and $K$ is a non-trivial knot in $Y$.
Then $g(K)+2d(Y)>0$ and the following bound holds:
\begin{equation}
\label{upbound}
p< \frac{4g(K)(g(K)+1)}{g(K)+2d(Y)}.
\end{equation}
\end{theorem}
We will prove the proof in Section~2.

By using the lower and upper bounds, we can understand that as $d(Y)$ is more increasing,
$g(K)$ must also become more increasing.
Note that lens spaces coming from L-space homology spheres are not empty other than $S^3$, for 
$L(22,3)$ is the $22$-Dehn surgery of Poincare homology sphere, which the correction term is $2$.
When $Y=S^3$, this theorem means Theorem~1 in \cite{[8]}.

In the case where an L-space homology sphere $Y$ is $d(Y)=2$, $Y$ can give
many lens spaces as demonstrated later on.
But when let $d(Y)\neq 0,2$,
there do not exist such lens spaces in the range of $p<1000$ and $|d(Y)|\le 40$ by Maple's computation.
For example $Y$ is $\Sigma(2,3,5)\#\Sigma(2,3,5)$.
It is not likely that lens space can be constructed from reducible manifolds.
We conjecture the following here.
\begin{conjecture}
Let $Y$ be any irreducible L-space homology sphere with $d(Y)\neq 0,2$.
Then there never exist knot $K\subset Y$ such that $Y_p(K)$ is a lens space for a positive integer $p$.
\end{conjecture}
The author proved this conjecture in the case where $Y$ is Poincar\'e homology sphere with reverse orientation,
nemely $d(Y)=-2$, (see \cite{[17]}).

In Section~4 we illustrate many lens spaces are obtained from L-space homology sphere with $d(Y)=2$.
In practice some of lens spaces are constructed from Poincar\'e homology sphere $\Sigma(2,3,5)$
along knots whose dual knots are $1$-bridge knots in the lens spaces.
We put a table of such lens spaces up to $p\le 2007$ in the last of the Section.
The L-space homology sphere with $d(Y)=2$ gives several infinite sequences of lens spaces except fewer 
examples.
In particular from Poincar\'e homology sphere infinite many lens spaces can be constructed
by Dehn surgeries.

The author does not know whether there exist lens spaces which are given rise to from both $S^3$ and $\Sigma(2,3,5)$.
Assume that exists such a lens space $L(p,q)$, the two knots $K,K'$ in $S^3,\Sigma(2,3,5)$ respectively.
Then $p$ is odd and the genus of $K'$ is $\frac{p+1}{2}$.
Moreover the Alexander polynomial of the two knots must have
the following relation $\Delta_{K'}(t)=\Delta_K(t)-(t^{\frac{p-1}{2}}+t^{-\frac{p-1}{2}})+(t^{\frac{p+1}{2}}+t^{-\frac{p+1}{2}})$,
by using Equations~(\ref{correction}) as explained later on.
\section{The exact triangle and the Alexander polynomial}
From this section on, we denote $Y$ as an L-space homology sphere.
If a positive $p$-Dehn surgery $Y_p(K)$ along knot $K$ is a lens space, then we have the short exact sequences
for any $0\neq i\in \Z/p\Z$,
\begin{equation}
0\to \mathop{\oplus}_{j\equiv i\mod p}HF^+(Y_0,j)\to HF^+(Y_p(K),Q(i))\to HF^+(Y)\to 0
\end{equation}
and for $i=0$,
\begin{equation}
0\to HF^+(Y)\to \mathop{\oplus}_{j\equiv 0\mod p}HF^+(Y_0,j)\to HF^+(Y_p(K),Q(i))\to 0.
\end{equation}
From this exact sequences for any $i$ the following formula is induced as in \cite{[3]}:
\begin{equation}
\label{correction}
d(Y)-d(Y_p(K),Q(i))+d(L(p,1),i)=2t_i(K),
\end{equation}
where $Q(i)=hi+c$, $h$ is the homology class of the dual knot $K^*$ of $K$ and 
$c=\frac{(h+1+p)(h-1)}{2}$ by the result in \cite{[6]}.
The each integer $t_i(K)$ is the $i$-th Turaev torsion of $Y_0(K)$.
It is positive when $Y_p(K)$ is L-space.

Here we use the identifications $\text{Spin}^c(L(p,q))\cong H_1(L(p,q))\cong H^2(L(p,q))\cong \Z/p\Z$.
By taking the summation of Equation~(\ref{correction}) over $i\in\Z/p\Z$.
we obtain 
\begin{eqnarray}
\label{Euler}
p\left(d(Y)+2\lambda(L(p,q))-2\lambda(L(p,1))\right)&=&2\sum_{i\in \Z}t_i(K)\nonumber\\
&=&2\sum_{i\in\Z}i^2a_i(K)\\
&=&\Delta_K''(t)|_{t=1},\nonumber
\end{eqnarray}
where the Alexander polynomial is normalized and symmetrized and Casson-Walker invariant $\lambda$ is
computed by the Rustamov's formula in \cite{[7]}:
$$\sum_{\frak{s}\in\text{Spin}^c(W)}\left(\chi(HF_{red}(W,\frak{s}))-
\frac{1}{2}d(W,\frak{s})\right)=|H_1(W,\Z)|\lambda(W),$$
where $W$ is any rational homology sphere.

By the same argument as the proof of Theorem~7.2 in \cite{[3]} the statement replaced $S^3$ with
$Y$ in that theorem also holds.

We can compute the coefficients of Alexander polynomial of the knot yielding
lens surgery as follows.
Let $[[\alpha,\beta]]$ be $[\alpha,\beta]\cap\Z$ and 
$[\gamma]_p$ be a reduction modulo $p$ of $\gamma$ with $0\le [\gamma]_p<p$.
We define $\Phi^k_{p,q}(h)$ as $\#\{j\in[[1,h']]|[qj-k]_p\in[[1,h]]\}$, where 
$h=[h]_p,h'=[h^{-1}]_p$.
For any class $i\in \Z/p\Z$ the reduced coefficient $\tilde{a}_i(K)$ of 
the Alexander polynomial of knot $K$ is
defined as $\sum_{j\equiv i\bmod p}a_j(K)$.
When $L(p,q)=Y_p(K)$, $\tilde{a}_{j}(K)$ is computed as follows:
\begin{equation}
\label{alex}
\tilde{a}_i(K)=-m+\Phi^{hi+c}_{p,q}(h).
\end{equation}
\section{The upper bound of the slope by the Seifert genus}
In this section we prove Theorem~\ref{main1}.
First, we introduce the two lemmas below.
\begin{lemma}
\label{degg}
Let $Y$ be an L-space homology sphere.
If $Y_p(K)$ is a lens space, the degree of the Alexander polynomial $\Delta_K(t)$ 
coincides with the Seifert genus $g(K)$.
\end{lemma}
\begin{lemma}
\label{os}
Let $Y$ be L-space homology sphere and $Y_p(K)$ be L-space.
Then $\Delta_K(t)$ has the form
$$\Delta_K(t)=(-1)^k+\sum_{j=1}^k(-1)^{k-j}(t^{n_j}+t^{-n_j})$$
for some increasing sequence of positive integers $0<n_1<n_2<\cdots<n_k$.
\end{lemma}
The proofs of these lemmas' are immedieately derived from an 
application to \cite{[9]} and \cite{[13]}.
Note that the property of Lemma~\ref{os} is also preserved for
the reduced polynomial $\displaystyle{\sum_{|i|<\frac{p}{2}}\tilde{a}_i(K)t^i}$
when $g(K)<\frac{p}{2}$ or 
$\displaystyle{\sum_{|i|<\frac{p}{2}}\tilde{a}_i(K)t^i+t^{\frac{p}{2}}+t^{-\frac{p}{2}}}$
when $g(K)=\frac{p}{2}$.

From this assumption $Y-K$ is irreducible and 
from the main theorem in \cite{[14]} and Lemma~\ref{os}, $K$ is a fiber knot.
This leads to the result of Lemma~\ref{degg}.

For the proof of Theorem~\ref{main1} we use the following proposition in \cite{[8]}.
Note that $\lambda$ is multiplied by $-1/2$ for the Casson-Walker invariant in \cite{[8]} and
the definition of orientations of lens spaces is the opposite to the one in \cite{[8]}.
\begin{proposition}[Proposition~2.4.\cite{[8]}]
\label{ras}
Suppose that $L$ is a lens space with $|H_1(L)|=p$, and that
$$2\lambda(L)-2\lambda(L(p,1))\le \frac{1}{4}(\frac{p}{4}-1).$$
Then $L$ is homeomorphic to one of $L(p,1)$, $L(p,2)$, or $L(p,3)$.
\end{proposition}
{\bf Proof of Theorem~\ref{main1})}
From the Fr\o yshov's inequality in \cite{[8]} and Equation~(\ref{Euler}), 
\begin{eqnarray}
2p(\lambda(L(p,q))-\lambda(L(p,1)))&=&2\sum_{i\ge 1}i^2a_i(K)-pd(Y)\nonumber\\
&\le&g(K)(g(K)+1)-pd(Y).\nonumber
\end{eqnarray}
We assume that
\begin{equation}
\label{ass}
p\ge \frac{4g(K)(g(K)+1)}{g(K)+2d(Y)}.
\end{equation}
Then we have $g(K)(g(K)+1)-pd(Y)\le \frac{p}{4}(\frac{p}{4}-1)$ by solving the quadratic inequality.
By Proposition~\ref{ras} $L(p,q)=L(p,1),L(p,2),$ or $L(p,3)$ hold.
L-space homology sphere $Y$ yielding such lens spaces must be $d(Y)=0$ or $2$, due to
Lemma~\ref{os} and the positivity of $t_i(K)$ for any $i$.

We consider the case of $Y_p(K)=L(p,1)$.
By Lemma~\ref{os} $d(Y)=0$ or $2$.
From $q=1$ we have $h=h'$ and we may assume that $h<\frac{p}{2}$ by replacing $h$ with $p-h$.
Hence Formula (\ref{alex}) says $\tilde{a}_i=-m$ for some integer $i$.
By Lemma~\ref{os} non-negative integer $m$ is $0$ or $1$.
When $m=0$, $h=h'=1$ nemely $\Delta_K(t)=1\bmod t^p-1$ holds.
$K$ is tirival from Lemma~\ref{degg} Therefore $Y=S^3$.
When $m=1$, $\tilde{a}_{-h'c}(K)=-1+h$ is $0$, $\pm 1$, or $2$ by applying Lemma~\ref{os}.
Thus we have $h=1,2,3$.\\
When $h=1$, $h^2=1$ hence $p=0$ and this is inconsistent.\\
When $h=2$, $h^2=4=p+1$ hence $p=3$. $\Delta_K(t)=1\bmod t^p-1$.\\
When $h=3$, $h^2=9=p+1$ hance $p=8$, $\Delta_K(t)=t^{-4}-t^{-3}+t-1+t-t^3+t^4$.\\
As a whole, the case of $L(p,q)=Y_p(K)$ classifies the following.
\begin{itemize}
\item $S^3_p(\text{trivial knot})=L(p,1)$ and $h=1$.
\item $Y_p(K_{1,p})=L(p,1)$, $d(Y)=2$, $p$ is odd, $h=1$\\
and $\Delta_{K_{1,p}}(t)=t^{-\frac{p+1}{2}}-t^{-\frac{p-1}{2}}+1-t^{\frac{p-1}{2}}+t^{\frac{p+1}{2}}$.
\item $Y_8(K_2)=L(8,1)$, $d(Y)=2$, $h=3$, $\Delta_{K_2}(t)=t^{-4}-t^{-3}+t^{-1}-1+t-t^3+t^4$.
\end{itemize}
In case of $Y_p(K)=L(p,2)$ or $L(p,3)$, analogous computation of $\tilde{a}_i$ and Lemma~\ref{degg}, \ref{os} classifies
the follows.
\begin{itemize}
\item $Y_7(K_3)=L(7,2)$, $d(Y)=0$ or $2$\\
and $\Delta_{K_3}(t)=t^{-1}-1+t$ or $t^{-4}-t^{-3}+t^{-1}-1+t-t^3+t^{4}$ respectively.
\item $Y_{11}(K_4)=L(11,3)$, $d(Y)=0$ or $2$, $h=5$
and $\Delta_{K_4}(t)=t^{-2}-t^{-1}+1-t+t^2$\\
or $t^{-6}-t^{-5}+t^{-2}-t^{-1}+1-t+t^2-t^5+t^{-6}$
respectively.
\item $Y_{13}(K_5)=L(13,3)$, $d(Y)=0$ or $2$, $h=5$
and $\Delta_{K_5}(t)=t^{-3}-t^{-2}+1-t^2+t^3$ or $t^{-7}-t^{-6}+t^{-3}-t^{-2}+1-t^2+t^3-t^{-6}+t^{-7}$
respectively.
\item $Y_{22}(K_6)=L(22,3)$ $d(Y)=2$, $h=5$\\
and $\Delta_{K_6}(t)=t^{-11}-t^{-10}+t^{-6}-t^{-5}+t^{-2}-1+t^2-t^5+t^6-t^{10}+t^{11}$.
\end{itemize}
Each of the cases does not satisfy Inequality~(\ref{ass}).
Therefore Inequality~(\ref{upbound}) is proven.\\
\hfill$\Box$
\section{A table of several lens surgeries over $Y$ with $d(Y)=2$}
In this section and the subsequent section $Y$ is an L-space homology sphere with $d(Y)=2$ for example
$Y$ is $\Sigma(2,3,5)$.
Since lens surgery by $K$ with $2g(K)-1=p$ includes the ones admitting lens surgery on $S^3$
we restrict our attention to lens surgery with $2g(K)-1<p$.
We may be able to construct a lens space form both integral Dehn surgeries of 
$S^3$ and $\Sigma(2,3,5)$ but ignore such a lens space here.
The homology class of the dual knot $K^*$ of $K$ is $h[l]\in H_1(L(p,q))$, where curve $l$ is a core loop of a
handlebody of genus one Heegaard splitting.
The set of classes $\mathcal{H}(p,K):=\{\pm h^{\pm}\}\subset \Z/p\Z$ is an invariant of the lens surgery which is independent of choices of two handlebodies of the
Heegaard decomposition and the orientation and we always consider any element in this set as the integers reduced to $\{1,2,\cdots,p\}$.

Any datum $(p,q,h,g)$ in Table~1, and 2 represent lens space $L(p,q)$ and 
the minimal representative $h$ in $\mathcal{H}(p,K)$ that the coefficient (\ref{alex}) 
satisfies Lemma~\ref{os}, and $t_i(K)$ is all non-negative, and 
$g$ is $\max\{i\in\{0,1,\cdots,\left[\frac{p}{2}\right]\}|-m+\Phi_{p,q}^{hi+c}(h)\neq 0\}$.


Let $L(p,q)=Y_p(K)$.
We say that $K$  is {\it 0-bridge knot} if $K$ is isotopic to a knot which 
lies on Heegaard surface of genus one 
Heegaard splitting, and {\it 1-bridge knot} if the knot is non-0-bridge knot and is the union of
two arcs embedded in the meridian disks of both handlebodies of genus one Heegaard splitting of the lens space.
This definition is based on \cite{[2]}.
From any triplicity $(p,q,h)$ we can uniquely find 0-bridge knot or 1-bridge knot having the datum.
\begin{theorem}
\label{main2}
Lens spaces $L(p,q)$ in Table~1 and 2 are constructed by $p$-Dehn surgery of knots in $\Sigma(2,3,5)$.
Moreover the dual knots are 1-bridge knots in $L(p,q)$.
\end{theorem}
Before the proof of this theorem we prove the following lemma.
\begin{lemma}
\label{lem}
Lens spaces in Table~1, 2 satisfy one of the list below for some $\ell\in\Z\setminus\{0\}$.
\begin{enumerate}
\item[\normalfont{a)}] $p=14\ell^2+7\ell+1$, $h=\pm(7\ell+2)^{\pm1} \mod p$, $2g=p+1-|\ell|$
\item[\normalfont{b)}] $p=20\ell^2+15\ell+3$, $h=\pm(5\ell+2)^{\pm1} \mod p$, $2g=p+1-|\ell|$
\item[\normalfont{c)}] $p=30\ell^2+9\ell+1$, $h=\pm(6\ell+1)^{\pm1} \mod p$, $2g=p+1-|\ell|$
\item[\normalfont{d)}] $p=42\ell^2+23\ell+3$, $h=\pm(7\ell+2)^{\pm1} \mod p$, $2g=p+1-|\ell|$
\item[\normalfont{d')}] $p=42\ell^2+47\ell+13$, $h=\pm(7\ell+4)^{\pm1} \mod p$, $2g=p+1-|\ell|$
\item[\normalfont{e)}] $p=52\ell^2+15\ell+1$, $h=\pm(13\ell+2)^{\pm1} \mod p$, $2g=p+1-|\ell|$
\item[\normalfont{e')}] $p=52\ell^2+63\ell+19$, $h=\pm(13\ell+8)^{\pm1} \mod p$, $2g=p+1-|\ell|$
\item[\normalfont{f)}] $p=54\ell^2+15\ell+1$, $h=\pm(27\ell+4)^{\pm1} \mod p$, $2g=p+1-|\ell|$
\item[\normalfont{f')}] $p=54\ell^2+39\ell+7$, $h=\pm(27\ell+10)^{\pm1} \mod p$, $2g=p+1-|\ell|$
\item[\normalfont{g)}] $p=69\ell^2+17\ell+1$, $h=\pm(23\ell+3)^{\pm1} \mod p$, $2g=p+1-2|\ell|$
\item[\normalfont{g')}] $p=69\ell^2+29\ell+3$, $h=\pm(23\ell+5)^{\pm1} \mod p$, $2g=p+1-2|\ell|$
\item[\normalfont{h)}] $p=85\ell^2+19\ell+1$, $h=\pm(17\ell+2)^{\pm1} \mod p$, $2g=p+1-2|\ell|$
\item[\normalfont{h')}] $p=85\ell^2+49\ell+7$, $h=\pm(17\ell+5)^{\pm1} \mod p$, $2g=p+1-2|\ell|$
\item[\normalfont{i)}] $p=99\ell^2+35\ell+3$, $h=\pm(11\ell+2)^{\pm1} \mod p$, $2g=p+1-2|\ell|$
\item[\normalfont{i')}] $p=99\ell^2+53\ell+7$, $h=\pm(11\ell+3)^{\pm1} \mod p$, $2g=p+1-2|\ell|$
\item[\normalfont{j)}] $p=120\ell^2+16\ell+1$, $h=\pm(12\ell+1)^{\pm1} \mod p$, $2g=p+1-2|\ell|$
\item[\normalfont{k)}] $p=120\ell^2+20\ell+1$, $h=\pm(20\ell+2)^{\pm1} \mod p$, $2g=p+1-2|\ell|$
\item[\normalfont{l)}] $p=120\ell^2+36\ell+3$, $h=\pm(12\ell+2)^{\pm1} \mod p$, $2g=p+1-2|\ell|$
\item[\normalfont{m)}] $p=120\ell^2+104\ell+22$, $h=\pm(12\ell+5)^{\pm1} \mod p$, $2g=p+1-|2\ell+1|$
\item[\normalfont{n)}] $L(191,34)$, $h=15$
\end{enumerate}
\end{lemma}
\proof
It is only necessary to prove that data $(p,q,h)$ given rise to by the 20 lists cover Table~1, 2 by direct computation.
The computation of genus of $K$ is due to Lemma~\ref{degg}.\\
\hfill$\Box$\\

\noindent
{\bf Proof of Theorem~\ref{main2})}
From each datum $(p,q,h)$ in Lemma~\ref{lem} we take the 1-bridge knot $K^*$ in $L(p,q)$.
The $-\tilde{a}_{-h'c-h'}(K)^*$-surgery in the sense of T.Saito in \cite{[15]} obtains a homology sphere $Y$.
Thus there exists a knot $K\subset Y$ such that we have $Y_p(K)=L(p,q)$.
The presentation of $\pi_1(Y)$ is
$$\left<x_1,x_2\big|\prod_{i=1}^{p}x_1x_2^{\delta_h(qi+1)},(\prod_{i=1}^{h'-1} x_1x_2^{\delta_h(qi+1)})x_1x_2^{-\tilde{a}_{h'c-h}}\right>$$
by \cite{[5]}, where $\delta_h:\Z/p\Z\to \{0,1\}$ is define to be $\delta_h(k)=1$ when $k=1,\cdots,$ or $h \mod p$, and $\delta_h(k)=0$ otherwise.

By deformating this group presentation directly we can easily give an isomorphism $\pi_1(Y)\cong \langle x,y|(xy)^2=x^3=y^5\rangle$.
The right hand side is the same as the fundamental group of $\Sigma(2,3,5)$.
Due to the celebrated resolution of Poincar\'e conjecture by G. Perelman in \cite{[16]}
, $Y$ is homeomorphic to $\Sigma(2,3,5)$.\\
\hfill$\Box$\\
The author would like to classify lens spaces obtained from Poincar\'e homology sphere.
But the author could not control type $n)$ in Lemma~\ref{lem}.
If you go on calculating more, then a new sequence may be discovered.

We here illustrate two figures, which are Figure~1 and 2.
Any point of Figure~1 represents a lens surgery over $\Sigma(2,3,5)$ with the slope $p\le 2007$.
The point $(h,p)$ means that there exists a knot such that $L(p,q)=\Sigma(2,3,5)_p(K)$ and
$h$ is the minimal in the set $\mathcal{H}(p,K)$ defined in Section~4.
Figure~2 represents hyperbolic lens surgeries over $S^3$ which are plotted in the order 
of slope from the smallest up to the same cardinarity as plots in Figure~1.
The horizontal and vertical axes mean the same as Figure~1.
Compared with Figure~1 and 2 it follows that the existence of lens surgeries 
over $\Sigma(2,3,5)$ are localized near four quadratic functions.
The right two of the them correspond to f), f'), g) and g').
To draw Figure~1, we referred to the last table of the article \cite{[2]}.
Here we give a rough conjecture on the basis of Figure~1.
\begin{table}[tbp]
\caption{Lens spaces with $p\le 711$ which homology sphere with $d(Y)=2$ yield.}\begin{center}
\small
\begin{tabular}{||r|r|r|r||r|r|r|r||r|r|r|r||} \hline
$p$ & $q$ & $h$ & $g$ & $p$ & $q$ & $h$ & $g$ & $p$ & $q$ & $h$ & $g$ \\ \hline
\phantom{00}8 & \phantom{0}1 & \phantom{0}3 & \phantom{00}4 & 221 & 127 & 41 & 109 & 442 & 157 & 77 & 220 \\ \hline
\phantom{0}22 & \phantom{0}3 & \phantom{0}5 & \phantom{0}11 & 228 & \phantom{0}61 & 17 & 113 & 445 & 186 & 39 & 221 \\ \hline
\phantom{0}38 & \phantom{0}7 & \phantom{0}7 & \phantom{0}19 & 239 & \phantom{0}67 & 28 & 119 & 445 & \phantom{0}84 & 23 & 221 \\ \hline
\phantom{0}40 & \phantom{0}9 & \phantom{0}7 & \phantom{0}20 & 243 & 133 & 43 & 120 & 449 & \phantom{0}80 & 23 & 223 \\ \hline
\phantom{0}43 & 15 & 12 & \phantom{0}21 & 244 & \phantom{0}45 & 17 & 121 & 450 & \phantom{0}79 & 23 & 224 \\ \hline
\phantom{0}53 & 11 & \phantom{0}8 & \phantom{0}26 & 246 & \phantom{0}43 & 17 & 122 & 463 & 211 & 40 & 229 \\ \hline
\phantom{0}67 & 14 & \phantom{0}9 & \phantom{0}33 & 247 & 134 & 58 & 123 & 469 & 107 & 24 & 233 \\ \hline
\phantom{0}68 & 13 & \phantom{0}9 & \phantom{0}34 & 249 & \phantom{0}94 & 29 & 123 & 497 & \phantom{0}79 & 24 & 247 \\ \hline
\phantom{0}70 & 11 & \phantom{0}9 & \phantom{0}35 & 250 & \phantom{0}39 & 17 & 124 & 509 & 116 & 25 & 253 \\ \hline
\phantom{0}71 & 38 & 16 & \phantom{0}35 & 253 & 141 & 30 & 125 & 513 & 112 & 25 & 255 \\ \hline
\phantom{0}87 & 13 & 10 & \phantom{0}43 & 263 & \phantom{0}61 & 18 & 130 & 514 & 139 & 41 & 256 \\ \hline
100 & 29 & 27 & \phantom{0}50 & 275 & \phantom{0}49 & 18 & 137 & 517 & 108 & 25 & 257 \\ \hline
101 & 21 & 18 & \phantom{0}50 & 294 & \phantom{0}67 & 19 & 146 & 521 & 201 & 42 & 259 \\ \hline
102 & 19 & 11 & \phantom{0}51 & 297 & \phantom{0}64 & 19 & 147 & 532 & \phantom{0}93 & 25 & 265 \\ \hline
103 & 18 & 11 & \phantom{0}51 & 298 & \phantom{0}13 & 19 & 148 & 532 & 309 & 85 & 265 \\ \hline
105 & 16 & 11 & \phantom{0}52 & 298 & \phantom{0}67 & 31 & 148 & 537 & 337 & 64 & 266 \\ \hline
106 & 37 & 19 & \phantom{0}52 & 301 & 176 & 64 & 150 & 547 & 295 & 44 & 271 \\ \hline
113 & 31 & 12 & \phantom{0}56 & 303 & 115 & 32 & 150 & 555 & 121 & 26 & 276 \\ \hline
125 & 19 & 12 & \phantom{0}62 & 311 & 168 & 49 & 154 & 571 & 202 & 66 & 283 \\ \hline
134 & 39 & 21 & \phantom{0}67 & 312 & \phantom{0}49 & 19 & 155 & 578 & 151 & 27 & 287 \\ \hline
137 & 30 & 13 & \phantom{0}68 & 316 & \phantom{0}65 & 33 & 156 & 583 & \phantom{0}93 & 26 & 290 \\ \hline
138 & 31 & 13 & \phantom{0}68 & 329 & \phantom{0}71 & 20 & 163 & 599 & 139 & 44 & 298 \\ \hline
139 & 30 & 13 & \phantom{0}69 & 337 & 188 & 51 & 167 & 610 & 351 & 91 & 304 \\ \hline
141 & 37 & 22 & \phantom{0}70 & 353 & \phantom{0}97 & 34 & 176 & 625 & 241 & 46 & 310 \\ \hline
145 & 51 & 44 & \phantom{0}72 & 376 & 145 & 71 & 187 & 633 & 151 & 28 & 314 \\ \hline
148 & 85 & 23 & \phantom{0}73 & 379 & 159 & 36 & 188 & 638 & \phantom{0}93 & 47 & 316 \\ \hline
159 & 37 & 14 & \phantom{0}79 & 383 & 101 & 22 & 190 & 673 & 473 & 72 & 334 \\ \hline
179 & 39 & 24 & \phantom{0}89 & 386 & 211 & 37 & 191 & 676 & 181 & 47 & 337 \\ \hline
187 & 69 & 50 & \phantom{0}93 & 411 & \phantom{0}73 & 22 & 204 & 706 & 135 & 29 & 351 \\ \hline
191 & 34 & 15 & \phantom{0}95 & 424 & 157 & 37 & 211 & 709 & 251 & 49 & 352 \\ \hline
197 & 51 & 26 & \phantom{0}97 & 428 & \phantom{0}89 & 23 & 212 & 710 & 131 & 29 & 353 \\ \hline
217 & 39 & 16 & 108 & 441 & 121 & 38 & 219 & 711 & 493 & 74 & 353 \\ \hline
\end{tabular}
\end{center}
\end{table} %

\begin{table}[tbp]
\caption{Lens spaces with $712\le p\le 2007$ which homology sphere with $d(Y)=2$ yield.}\begin{center}
\small
\begin{tabular}{||r|r|r|r||r|r|r|r||r|r|r|r||} \hline
$p$ & $q$ & $h$ & $g$ & $p$ & $q$ & $h$ & $g$ & $p$ & $q$ & $h$ & $g$ \\ \hline
\phantom{0}715 & 199 & \phantom{0}98 & 356 & 1103 & 291 & \phantom{0}60 & 550 & 1552 & \phantom{0}849 & 145 & 774 \\ \hline
\phantom{0}736 & 393 & \phantom{0}51 & 365 & 1129 & 240 & \phantom{0}37 & 562 & 1563 & \phantom{0}640 & \phantom{0}73 & 778 \\ \hline
\phantom{0}739 & 161 & \phantom{0}30 & 367 & 1135 & 234 & \phantom{0}37 & 565 & 1583 & \phantom{0}334 & 110 & 787 \\ \hline
\phantom{0}767 & 133 & \phantom{0}30 & 382 & 1141 & 421 & \phantom{0}62 & 568 & 1618 & \phantom{0}149 & \phantom{0}75 & 804 \\ \hline
\phantom{0}773 & 181 & \phantom{0}50 & 385 & 1162 & 253 & 125 & 579 & 1634 & \phantom{0}427 & \phantom{0}73 & 815 \\ \hline
\phantom{0}789 & 172 & \phantom{0}31 & 392 & 1163 & 149 & \phantom{0}38 & 578 & 1641 & \phantom{0}340 & 112 & 816 \\ \hline
\phantom{0}790 & 171 & \phantom{0}31 & 393 & 1168 & 201 & \phantom{0}37 & 582 & 1653 & \phantom{0}283 & \phantom{0}44 & 824 \\ \hline
\phantom{0}796 & 165 & \phantom{0}31 & 396 & 1171 & 321 & \phantom{0}63 & 582 & 1717 & \phantom{0}307 & 152 & 856 \\ \hline
\phantom{0}805 & 211 & 104 & 401 & 1173 & 814 & \phantom{0}95 & 583 & 1727 & \phantom{0}389 & \phantom{0}46 & 860 \\ \hline
\phantom{0}813 & 211 & \phantom{0}32 & 404 & 1191 & 253 & \phantom{0}38 & 593 & 1742 & \phantom{0}283 & \phantom{0}45 & 868 \\ \hline
\phantom{0}823 & 340 & \phantom{0}53 & 409 & 1198 & 631 & \phantom{0}65 & 595 & 1758 & \phantom{0}451 & \phantom{0}47 & 875 \\ \hline
\phantom{0}828 & 133 & \phantom{0}31 & 412 & 1223 & 848 & \phantom{0}97 & 608 & 1772 & \phantom{0}925 & \phantom{0}79 & 881 \\ \hline
\phantom{0}841 & 107 & \phantom{0}54 & 417 & 1226 & 257 & \phantom{0}63 & 611 & 1779 & \phantom{0}337 & \phantom{0}46 & 886 \\ \hline
\phantom{0}873 & 151 & \phantom{0}32 & 435 & 1243 & 201 & \phantom{0}38 & 619 & 1783 & \phantom{0}427 & \phantom{0}76 & 889 \\ \hline
\phantom{0}878 & 129 & \phantom{0}33 & 436 & 1276 & 291 & 131 & 636 & 1803 & \phantom{0}406 & \phantom{0}47 & 898 \\ \hline
\phantom{0}893 & 237 & \phantom{0}54 & 445 & 1285 & 336 & \phantom{0}66 & 639 & 1807 & \phantom{0}309 & \phantom{0}46 & 901 \\ \hline
\phantom{0}919 & 379 & \phantom{0}56 & 457 & 1298 & 223 & \phantom{0}39 & 647 & 1811 & 1247 & 118 & 901 \\ \hline
\phantom{0}925 & 519 & 112 & 461 & 1331 & 135 & \phantom{0}68 & 661 & 1841 & \phantom{0}561 & \phantom{0}78 & 917 \\ \hline
\phantom{0}938 & 151 & \phantom{0}33 & 467 & 1376 & 361 & \phantom{0}67 & 686 & 1849 & \phantom{0}360 & \phantom{0}47 & 921 \\ \hline
\phantom{0}953 & 505 & \phantom{0}58 & 473 & 1377 & 223 & \phantom{0}40 & 686 & 1853 & \phantom{0}189 & \phantom{0}48 & 922 \\ \hline
\phantom{0}975 & 181 & \phantom{0}34 & 485 & 1379 & 302 & \phantom{0}41 & 686 & 1855 & \phantom{0}319 & 158 & 925 \\ \hline
\phantom{0}991 & 265 & \phantom{0}87 & 492 & 1403 & 361 & \phantom{0}42 & 698 & 1857 & \phantom{0}352 & \phantom{0}47 & 925 \\ \hline
\phantom{0}999 & 226 & \phantom{0}35 & 497 & 1408 & 273 & \phantom{0}41 & 701 & 1873 & 1289 & 120 & 932 \\ \hline
1004 & 233 & \phantom{0}57 & 500 & 1414 & 267 & \phantom{0}41 & 704 & 1887 & \phantom{0}406 & \phantom{0}80 & 939 \\ \hline
1021 & 301 & \phantom{0}58 & 508 & 1426 & 783 & 139 & 711 & 1900 & \phantom{0}289 & \phantom{0}47 & 947 \\ \hline
1027 & 189 & \phantom{0}35 & 511 & 1437 & 589 & \phantom{0}70 & 715 & 1933 & \phantom{0}163 & \phantom{0}82 & 961 \\ \hline
1027 & 573 & 118 & 512 & 1447 & 317 & \phantom{0}42 & 720 & 1963 & \phantom{0}511 & \phantom{0}80 & 979 \\ \hline
1033 & 192 & \phantom{0}35 & 514 & 1471 & 771 & \phantom{0}72 & 731 & 1985 & \phantom{0}416 & \phantom{0}49 & 989 \\ \hline
1037 & 271 & \phantom{0}89 & 515 & 1488 & 169 & \phantom{0}43 & 740 & 1993 & \phantom{0}408 & \phantom{0}49 & 993 \\ \hline
1057 & 239 & \phantom{0}36 & 526 & 1513 & 285 & \phantom{0}43 & 760 & 2001 & \phantom{0}721 & \phantom{0}82 & 997 \\ \hline
1072 & 121 & \phantom{0}61 & 532 & 1526 & 323 & \phantom{0}43 & 760 &  &  &  &  \\ \hline
1088 & 281 & \phantom{0}37 & 541 & 1534 & 315 & \phantom{0}43 & 764 &  &  &  &  \\ \hline
\end{tabular}
\end{center}
\end{table} %

\begin{conjecture}
Suppose $L(p,q)=\Sigma(2,3,5)_p(K)$.
For $h \in\mathcal{H}(p,K)$, either of the following 6 patterns holds:
\begin{enumerate}
\item $L(p,q)=L(54\ell^2+15\ell+1,27\ell^2+21\ell+3)$ for $\ell\in \Z\setminus \{0\}$,
\item $L(p,q)=L(54\ell^2+39\ell+7,27\ell^2+33\ell+9)$ for $\ell\in \Z\setminus \{0\}$,
\item $L(p,q)=L(69\ell^2+17\ell+1,46\ell^2+19\ell+2)$ for $\ell\in \Z\setminus \{0\}$,
\item $L(p,q)=L(69\ell^2+29\ell+3,46\ell^2+27\ell+4)$ for $\ell\in \Z\setminus \{0\}$,
\item $3.21\le \frac{h^2}{p}\le3.61$,
\item $1.15\le \frac{h^2}{p}\le1.28$.
\end{enumerate}
\end{conjecture}
\section*{Acknowledgments}
The author would like express the gratitude for Professor Masaaki Ue and 
would thank for T. Kadokami's organizing seminars and useful communication there with 
N. Maruyama and T .Sakai as well.
The author is grateful for J. Rasmussen's pointing out mistakes of the list in Lemma~\ref{lem}.
The author was supported by the Grant-in-Aid for JSPS Fellows for Young 
Scientists(17-1604).

\begin{figure}[h]
\begin{center}
\includegraphics{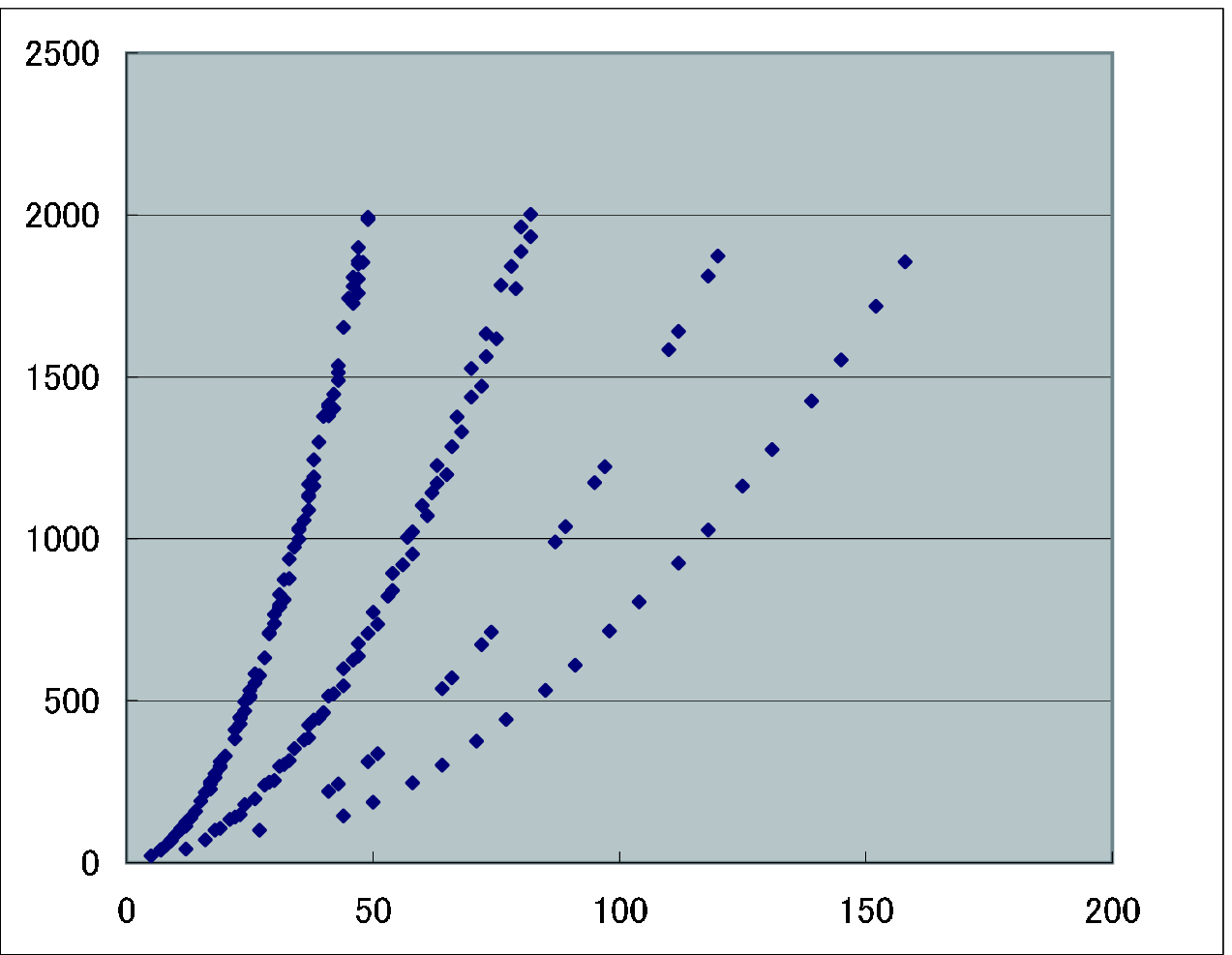}
\caption{$h-p$ graph in $\Sigma(2,3,5)$ case.}
\end{center}
\label{graph}
\begin{center}
\includegraphics{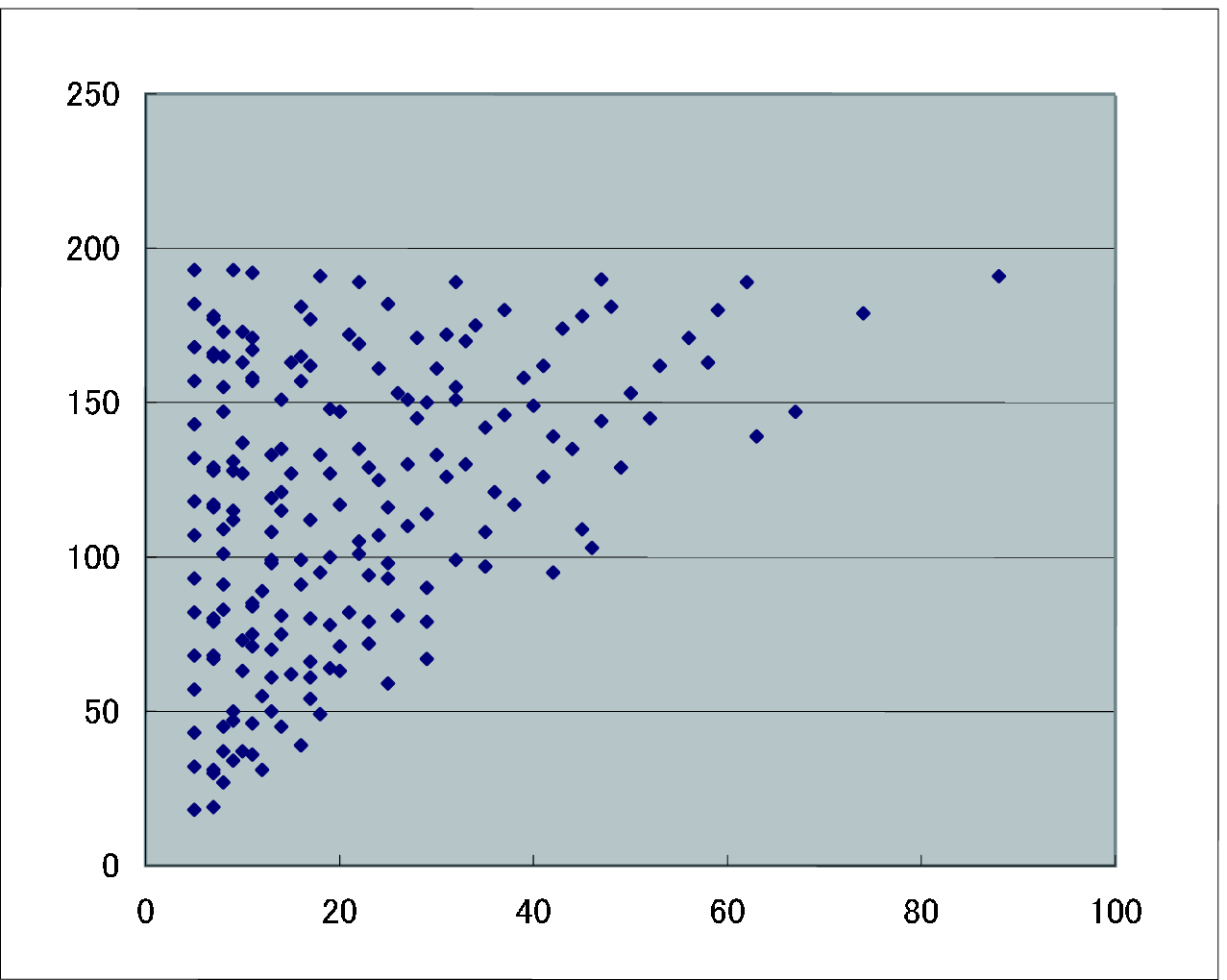}
\caption{$h-p$ graph in $S^3$ case.}
\label{S3graph}
\end{center}
\end{figure}


\begin{thebibliography}{000}
\bibitem{[2]} J. Berge, {\it Some knots with surgeries yielding lens spaces}, unpublished manuscript.

\bibitem{[1]} R. Fintushel and R. Stern, {\it Constructing lens spaces by surgery on knots},
Math. Z. Vol. 175, no. 1 February, (1980), 33-51

\bibitem{[4]} H. Goda and M. Teragaito, {\it Dehn surgeries on knots which yield lens spaces and genera of knots},
Math. Proc. Cambridge Philos. Soc. 129 (2000), 501-515.

\bibitem{[10]} T. Kadokami, {\it On the Alexander polynomial satisfying Ozsv\'ath Szab\'o's condition for lens surgery},
preprint

\bibitem{[12]} T. Kadokami, and Y. Yamada, {\it A deformation of the Alexander polynomials of knots yielding lens spaces},
Bull. of Austral. Math. Soc.

\bibitem{[9]} P. Kronheimer, T. Mrowka, P. Ozsv\'ath, and Z. Szab\'o, {\it Monopoles and lens space surgeries}, arXiv:math.GT/0310164

\bibitem{[11]}  K. Ichihara, T. Saito, and M. Teragaito, {\it Alexander polynomials of doubly primitive knots}, 
Proc. Amer. Math. Soc. 135 (2007), 605-615 

\bibitem{[3]} P. Ozsv\'ath and Z. Szab\'o, {\it Absolutely graded Floer homologies and intersection forms for four-manifolds with boundary}, 
Adv. Math. 173 (2003), no. 2, 179--261.

\bibitem{[9]} P. Ozsv\'ath and Z. Szab\'o, {\it Holomorphic disks and genus bounds},
Geometry \& Topology 8(2004) 311--334.

\bibitem{[13]} P. Ozsv\'ath and Z. Szab\'o, {\it On knot Floer homology and lens surgery},
Topology 
Volume 44, Issue 6 , November 2005, Pages 1281-1300 

\bibitem{[16]} G. Perelman, {\it  The entropy formula for the Ricci flow and its geometric applications }
 arXiv:math/0211159

\bibitem{[8]} J. Rasmussen, {\it Lens space surgeries and a conjecture of Goda and Teragaito}
Geometry \& Topology Vol.8(2004) 1013--1031

\bibitem{[7]} R. Rustamov, {\it Surgery formula for the renormalized Euler characteristic of Heegaard Floer homology},
arXiv:math.GT/0409294

\bibitem{[15]} T. Saito, {\it Dehn surgery and (1,1)-knots in lens spaces},
Topology Appl. 154 (2007), no. 7, 1502-1515

\bibitem{[17]} M. Tange, {\it  On the non-existence of lens space surgery structure},
arXiv:0707.0197
 
\bibitem{[6]} M. Tange, {\it Ozsv\'ath-Szab\'o's correction term of lens surgery},
preprint.
\bibitem{[5]} M. Tange, {\it The dual knot of lens surgery and Alexander polynomial},
preprint.

\bibitem{[14]} Y, Ni, {\it Knot Floer homology detects fiberd knots},
arXiv:math.GT/0607156
\end{thebibliography}
\end{document}